\documentclass[onecolumn,12pt,a4paper,draftcls]{IEEEtran}
\ifCLASSINFOpdf
\else
\fi

\usepackage[cmex10]{amsmath}

\usepackage[tight,footnotesize]{subfigure}
\usepackage{graphicx,amsfonts,amssymb,pstricks}
    \newcommand{\Mat}{\mathbf}
    \newcommand{\Cmplx}{\mathbb{C}}
    
    \newcommand{\Set}{\mathcal}
    \newcommand{\Exptn}[1]{\mathsf{E}\left\{#1\right\}}
	\newcommand{\Tr}[1]{\mathsf{Tr}\left[#1\right]}
	\newcommand{\Min}[1]{\mathop{\mbox{Minimize}}_{#1}}	
	\newcommand{\STo}{\mbox{Subject to}}
    \DeclareMathAlphabet{\mathitbf}{OML}{cmm}{b}{it}
    \DeclareMathSymbol{\alpha}{\mathalpha}{letters}{"0B}
    \DeclareMathSymbol{\beta}{\mathalpha}{letters}{"0C}
    \DeclareMathSymbol{\gamma}{\mathalpha}{letters}{"0D}
    \DeclareMathSymbol{\delta}{\mathalpha}{letters}{"0E}
    \DeclareMathSymbol{\epsilon}{\mathalpha}{letters}{"0F}
    \DeclareMathSymbol{\zeta}{\mathalpha}{letters}{"10}
    \DeclareMathSymbol{\eta}{\mathalpha}{letters}{"11}
    \DeclareMathSymbol{\theta}{\mathalpha}{letters}{"12}
    \DeclareMathSymbol{\iota}{\mathalpha}{letters}{"13}
    \DeclareMathSymbol{\kappa}{\mathalpha}{letters}{"14}
    \DeclareMathSymbol{\lambda}{\mathalpha}{letters}{"15}
    \DeclareMathSymbol{\mu}{\mathalpha}{letters}{"16}
    \DeclareMathSymbol{\nu}{\mathalpha}{letters}{"17}
    \DeclareMathSymbol{\xi}{\mathalpha}{letters}{"18}
    \DeclareMathSymbol{\pi}{\mathalpha}{letters}{"19}
    \DeclareMathSymbol{\rho}{\mathalpha}{letters}{"1A}
    \DeclareMathSymbol{\sigma}{\mathalpha}{letters}{"1B}
    \DeclareMathSymbol{\tau}{\mathalpha}{letters}{"1C}
    \DeclareMathSymbol{\upsilon}{\mathalpha}{letters}{"1D}
    \DeclareMathSymbol{\phi}{\mathalpha}{letters}{"1E}
    \DeclareMathSymbol{\chi}{\mathalpha}{letters}{"1F}
    \DeclareMathSymbol{\psi}{\mathalpha}{letters}{"20}
    \DeclareMathSymbol{\omega}{\mathalpha}{letters}{"21}
    \DeclareMathSymbol{\varepsilon}{\mathalpha}{letters}{"22}
    \DeclareMathSymbol{\vartheta}{\mathalpha}{letters}{"23}
    \DeclareMathSymbol{\varpi}{\mathalpha}{letters}{"24}
    \DeclareMathSymbol{\varrho}{\mathalpha}{letters}{"25}
    \DeclareMathSymbol{\varsigma}{\mathalpha}{letters}{"26}
    \DeclareMathSymbol{\varphi}{\mathalpha}{letters}{"27}
    \DeclareMathSymbol{\Gamma}{\mathalpha}{letters}{"00}
    \DeclareMathSymbol{\Delta}{\mathalpha}{letters}{"01}
    \DeclareMathSymbol{\Theta}{\mathalpha}{letters}{"02}
    \DeclareMathSymbol{\Lambda}{\mathalpha}{letters}{"03}
    \DeclareMathSymbol{\Xi}{\mathalpha}{letters}{"04}
    \DeclareMathSymbol{\Pi}{\mathalpha}{letters}{"05}
    \DeclareMathSymbol{\Sigma}{\mathalpha}{letters}{"06}
    \DeclareMathSymbol{\Upsilon}{\mathalpha}{letters}{"07}
    \DeclareMathSymbol{\Phi}{\mathalpha}{letters}{"08}
    \DeclareMathSymbol{\Psi}{\mathalpha}{letters}{"09}
    \DeclareMathSymbol{\Omega}{\mathalpha}{letters}{"0A}

\begin{document}
\title{Robust Downlink Beamforming in Multiuser MISO Cognitive Radio Networks}

% use a multiple column layout for up to three different affiliations
%\author{
%    \IEEEauthorblockN{Ebrahim A. Gharavol}
%        \IEEEauthorblockA{
%            Department of Electrical and Computer Engineering\\
%            National University of Singapore\\
%            21 Lower Kent Ridge Road, Singapore 119077 \\
%            Email: {\tt Ebrahim@nus.edu.sg}
%        }
%    \and
%    \IEEEauthorblockN{Ying-Chang Liang}
%        \IEEEauthorblockA{
%            Institute of Infocomm Research,\\
%            Agency for Science, Technology and Research (A*STAR),\\
%            1 Fusionpolis Way, Singapore 138632\\
%            Email: {\tt ycliang@i2r.a-star.edu.sg}
%        }
%} % '\author'


\author{
   \IEEEauthorblockN{
       Ebrahim A. Gharavol\IEEEauthorrefmark{1},
       Ying-Chang Liang\IEEEauthorrefmark{2}, and
       Koenraad Mouthaan\IEEEauthorrefmark{1}
   }

   \IEEEauthorblockA{
       \IEEEauthorrefmark{1}
       		Department of Electrical and Computer Engineering,
            National University of Singapore,
            21 Lower Kent Ridge Road, Singapore 119077 \\
            Email: {\tt \{Ebrahim,k.mouthaan\}@nus.edu.sg}\\
   }
   \IEEEauthorblockA{
       \IEEEauthorrefmark{2}
       		Institute of Infocomm Research, Agency for Science, Technology and Research (A*STAR),
            1 Fusionpolis Way, Singapore 138632\\
            Email: {\tt ycliang@i2r.a-star.edu.sg}\\
   }
} %'\author'
%\journal{IEEE Transactions on Wireless Communications (Submitted)}


\maketitle
%\doublespacing

\begin{abstract}
This paper studies the problem of robust downlink beamforming design in a multiuser Multi-Input Single-Output (MISO) Cognitive Radio Network (CR-Net) in which multiple Primary Users (PUs) coexist with multiple Secondary Users (SUs). Unlike conventional designs in CR-Nets, in this paper it is assumed that the Channel State Information (CSI) for all relevant channels is imperfectly known, and the imperfectness of the CSI is modeled using an Euclidean ball-shaped uncertainty set. Our design objective is to minimize the transmit power of the SU-Transmitter (SU-Tx) while simultaneously targeting a lower bound on the received Signal-to-Interference-plus-Noise-Ratio (SINR) for the SU's, and imposing an upper limit on the Interference-Power (IP) at the PUs. The design parameters at the SU-Tx are the beamforming weights, i.e. the precoder matrix.
The proposed methodology is based on a worst case design scenario through which the performance metrics of the design are immune to variations in the channels. We propose three approaches based on convex programming for which efficient numerical solutions exist. Finally, simulation results are provided to validate the robustness of the proposed methods.
\end{abstract}
%\IEEEpeerreviewmaketitle

\section{Introduction}

A Cognitive Radio Network (CR-Net) \cite{mitola1999}, \cite{haykin2005} is an intelligent solution to the spectrum scarcity problem.
In a CR-Net, the Secondary Users (SUs) are allowed to operate within the service range of the Primary Users (PUs), though the PUs have higher priority in utilizing the spectrum. There are two types of CR-Nets: opportunistic CR-Nets for which the SUs sense the spectrum and try to utilize the unused channels when they are not occupied by PUs; and concurrent CR-Nets in which SUs are allowed to use the spectrum even when PUs are active, provided that the amount of interference power to each PU is kept below a certain threshold  \cite{haykin2005}.
Hereafter, we will focus on concurrent CR-Nets.

Fig.~\ref{fig:CRNET}-a illustrates the downlink scenario of multiuser multiple-input single-output (MISO) CR-Net with $K$ SUs coexist with $L$ PUs. There are two sets of relevant channel state information (CSI) which play important roles in the system design: one set describes the channels between SU-Transmitter (SU-Tx) and SU-Receivers (SU-Rx's) while the other set describes the channels between SU-Tx and PU-Receivers (PU-Rx's). For simplicity, we term the first set of CSI as SU-link CSI and the second set as PU-link CSI. When PUs are inactive, the system becomes conventional multiuser MISO system, and SU-link CSI is needed for transmission design. This knowledge is usually acquired through  transmitting pilot symbols from SU-Tx to SU-RXs, and feeding back the estimated CSI from SU-Rxs to SU-Tx. In practice, however, because of the time variant nature of wireless channels, it is not possible to acquire the CSI perfectly, either due to channel estimation error and/or feedback error. On the other hand, when PUs are active, PU-link CSI is further needed at SU-Tx for the purpose of controlling interferences at the PU-Rx's. This CSI knowledge has to be acquired by SU-Tx through environmental learning \cite{gao-ICASSP2009}, which again will introduce errors in CSI. In this paper we consider the transmit design for a multiuser MISO CR-Net with uncertain CSI in both SU-link and PU-link.

Previously in conventional radio network design, ad-hoc methods, such as diagonal loading \cite{bertsekas1995}, were exploited in the design procedure of robust beamforming systems.
Quite recently these designs are based on well-reputed mathematical methodologies, such as the systematical worst case designs \cite{vorobyov2003}-\cite{lorenz2005}.
These methods deal with a Minimum Variance Distortionless Response (MVDR) problem in the signal processing domain and show that the problem may be recast as a Second-Order Cone Program (SOCP) \cite{boyd2004}.
Also, it was shown that this worst case design scenario is equivalent to an adaptive diagonal loading \cite{vorobyov2003}.
One of the first worst case designs was published by Bengtsson and Otterstten \cite{bengtsson1999}-\cite{bengtsson2001}.
They showed that the robust maximization of SINR would lead to a Semi-Definite Program (SDP) \cite{boyd2004}, after a simple Semidefinite Relaxation (SDR).
Sharma, {\em et al.,} \cite{sharma2008} developed a model to cover the Positive Semi-definiteness (PSD) of the  channel covariance matrix.
They proposed two SDPs, a conventional SDP and a SDP based on an iterative algorithm.
Also Mutapcic {\em et al.,} \cite{mutapcic2007} proposed a new tractable method to solve the robust downlink beamforming.
Their method is based on the cutting set algorithm which is also an iterative method.
Also \cite{vucic2009}-\cite{botros2009} are targeting the robust design of a beamforming system using the worst case scenarios for Quality-of-Service (QoS) constraints.

Quite a few works are published on the robust design for CR-Nets \cite{zhang2008}, \cite{zhang2009}, \cite{zhi2008} and \cite{cumanan2008}.
Zhang {\em et al.} \cite{zhang2008}, \cite{zhang2009} have studied such a CR-Net from an information theoretic perspective. The CR-Net considered in \cite{zhang2008}, \cite{zhang2009} consists of one PU-Rx and one SU-Rx, and the SU-link CSI is assumed to be perfectly known, but the PU-link CSI has uncertainty. A duality theory was developed to cope with the CSI imperfectness. Additionally, the authors proposed an analytic solution for this case.
Also, Zhi {\em et al.} \cite{zhi2008} designed a robust beamformer for a CR-Net, where the system setup is the same as in \cite{zhang2009}, however there may be some uncertainty in both the channel covariance matrix as well as the antenna manifold. Finally, Cumanan {\em et al.} \cite{cumanan2008} considered a CR-Net having multiple PUs and only one SU. In this work, both channels are assumed to be imperfect. They also used the worst case design method to come up with a convex problem that can be solved efficiently.

In this paper, we consider a downlink system of a CR-Net with multiple PUs and multiple SUs whose relevant CSI is imperfectly known. The imperfectness of the CSI is modeled using an Euclidean ball. Our design objective is to minimize the transmit power of the SU-Tx while simultaneously targeting a lower bound on the received Signal-to-Interference-plus-Noise-Ratio (SINR) for the SUs, and imposing an upper limit on the Interference-Power (IP) at the PUs. The design parameters at the SU-Tx are the beamforming weights, i.e. the precoder matrix.
The proposed methodology is based on a worst case design scenario through which the performance metrics of the design are immune to variations in the channels. We propose three approaches based on convex programming for which efficient numerical solutions exist. In the first approach, the worst case SINR is derived through using loose upper and lower bounds on the terms appearing in the numerator and denominator of the SINR. Working in this line, SDP is developed which provides us the robust beamforming coefficients. In the second approach, the minimum SINR is found through minimizing its numerator while maximizing its denominator. Different from the first method, we chose exact upper and lower bounds on the previously mentioned terms. This approach does not lead to a SDP, but the resulting problem is still convex and may be solved efficiently.
Finally, in our third approach, we find the exact minimum of SINR directly, and this method is also a general convex optimization problem.

The rest of the paper is organized as follows.
In Section \ref{sysModel}, the model of our studied system is described.
In Section \ref{probFormulation}, the robust design of a multiuser MISO CR-Net with multiple SUs and multiple PUs is considered.
In Section \ref{method1}, we show that the resulting optimization problem, using loose upper and lower bounds, is a SDP. In Sections \ref{method2} and \ref{method3} we propose two more general problem formulations based on stricter bound and exact bound on the minimum value of SINR, respectively.
In Section \ref{simResults} the simulation results that demonstrate the robustness of the proposed schemes are presented.
Finally, in Section \ref{conclusion} we conclude the paper.

{\bf Notations:}
Matrices and vectors are typefaced using slanted bold uppercase and lowercase letters, respectively.
Conjugate and conjugate transpose of the matrix $\Mat{A}$ are denoted as the $\Mat{A}^\dag$ and $\Mat{A}^*$,  respectively.
The trace of a matrix is annotated using $\Tr{\cdot}$.
Positive semi-definiteness of the matrix $\Mat{A}$ is depicted using $\Mat{A} \succeq 0$.
The symbol ``$\triangleq$'' means ``defined as''.
$\Cmplx^{m\times n}$ is used to describe the complex space of $m \times n$ matrices.
A Zero-Mean Circularly Symmetric Complex Gaussian (ZMCSCG) random variable with the variance of $\sigma^2$ is denoted using $\Set{CN}(0,\sigma^2)$.
For a vector like $\Mat{x}$, $\|\Mat{x}\|$ is the Euclidean norm while the norm of a matrix like $\|\Mat{A}\|$ is the Frobenius norm.
To show the differentiation of a function, $f$, with respect to some of its parameters, $\Mat{a}$, $\nabla_{\Mat{a}} f(\cdot)$ is used.
Finally, mathematical expectation is described as $\Exptn{\cdot}$.

\section{System Model\label{sysModel}}

\begin{figure}[p]
	\centering
    \begin{tabular}{c}
    	\includegraphics[width=0.70\textwidth]{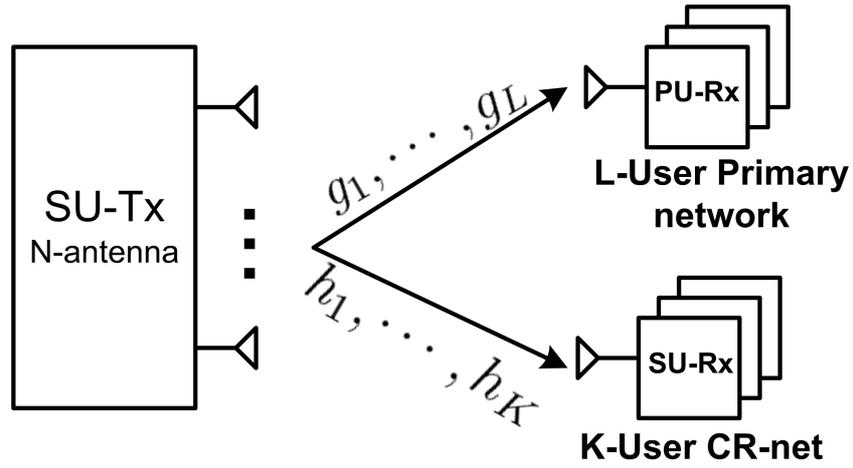} 	\\
    	\scriptsize (a) Symbolic Representation							\\
    	\includegraphics[width=0.90\textwidth]{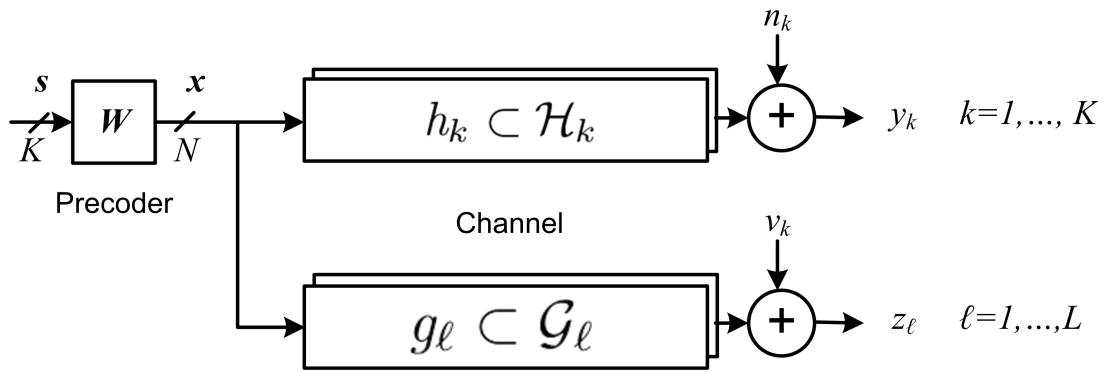}	\\
    	\scriptsize (b) Signal Processing Model
	\end{tabular}
    \caption{A Typical Multiuser MISO CR-Net System with Uncertain CSI}
    \label{fig:CRNET}
\end{figure}

Fig.~\ref{fig:CRNET} shows the downlink scenario of a multiuser MISO CR-Net coexisting with a Primary-Radio Network (PR-Net) having $L$ PUs each equipped with a single antenna. The SU-Tx equipped with $N$ antennas transmits independent symbols, $s_k$, to $K$ different single antenna SUs, $\{s_k \in \Cmplx \}_{k=1}^{K}$.
It is assumed that the transmitted symbols are all Gaussian, zero-mean and independent, i.e., $s_k \sim \Set{CN}(0,1)$.
Each symbol is precoded by a weight vector, $\{\Mat{w}_k \in \Cmplx^{N\times 1}\}_{k=1}^{K}$, resulting in a vector signal, $\{ \Mat{s}_k = \Mat{w}_k s_k \}_{k=1}^{K}$, for each one. It is known that $\Mat{s}_k \sim \Set{CN}(\Mat{0},\Mat{Q}_k)$, where $\Mat{Q}_k$ is the covariance matrix of $\Mat{s}_k$; $\Mat{Q}_k = \mathsf{E}\{\Mat{s}_k\Mat{s}_k^\dag\}=\Mat{w}_k\Mat{w}_k^\dag \succeq 0$ and $\Mat{0}$ is the zero vector.
The channel from SU-Tx to each SU-Rx is determined using a complex-valued vector, $\{ \Mat{h}_k \in \Cmplx^{N\times1} \}_{k=1}^{K}$ which is not perfectly known and there is some kind of uncertainty in channel gains.
This uncertainty is described using an uncertainty set, $\Set{H}_k$ which is defined as an Euclidean ball
\begin{equation}
	\Set{H}_k = \{ \Mat{h} | \| \Mat{h}-\Mat{\tilde{h}}_k \| \leq \delta_k\}.
\end{equation}
In this definition, the ball is centered around the actual value of the channel vector, $\Mat{\tilde{h}}_k$, and the radius of the ball is determined by $\delta_k$ which is a positive constant.
Using this notion, the channel is modeled as
\begin{equation}
	\Mat{h}_k = \Mat{\tilde{h}}_k + \Mat{a}_k,
\end{equation}	
where $\Mat{a}_k$ is a norm-bounded uncertainty vector, $\| \Mat{a}_k \| \leq \delta_k$.

The SU-Tx combines the signals and transmits the combination, $\Mat{x}$,
\begin{equation}
	\Mat{x} = \sum_{k=1}^{K} \Mat{s}_k = \Mat{Ws},
\label{eq:}
\end{equation}
where $\Mat{s} = [s_1,\cdots, s_K]^T \in \Cmplx^{K\times1}$ contains the transmitted symbols and as we know, $\Mat{s} \sim \Set{CN}(\Mat{0},\Mat{I}_K)$, also $\Mat{W} = [\Mat{w}_1, \cdots, \Mat{w}_K] \in \Cmplx^{N\times K}$, is called the precoding matrix.
For $\Mat{x}$ we know that $\Mat{x} \sim \Set{CN}(\Mat{0},\Mat{Q})$, where $\Mat{Q} = \mathsf{E}\{ \Mat{x}\Mat{x}^\dag \} = \Mat{WW}^\dag \succeq 0$.
The design objective is to determine this precoding matrix $\Mat{W}$ based on certain criteria that we will discuss in the next few paragraphs.

The channel from the SU-Tx to a PU-Rx is also defined using a complex valued vector, i.e., $\{ \Mat{g}_\ell \in \Cmplx^{N\times1} \}_{\ell=1}^{L}$.
Here it is assumed that the CSI for these users is also uncertain.
We use the same notation to describe the uncertainty for these channels.
The uncertainty is defined using a set, $\Set{G}_\ell$, which is
\begin{equation}
	\Set{G}_\ell=\{ \Mat{g}| \|\Mat{g} - \Mat{\tilde{g}}_\ell \| \leq \eta_\ell\}.
\end{equation}
Equivalently, we may write
\begin{equation}
	\Mat{g}_\ell = \Mat{\tilde{g}}_\ell + \Mat{b}_\ell,
\end{equation}
where $\Mat{b}_\ell$ is a norm-bonded uncertain vector, $\|\Mat{b}_\ell\| \leq \eta_\ell$ and $\Mat{\tilde{g}}_\ell$ is the actual value of the channel.

\section{Problem Formulation\label{probFormulation}}
For the system depicted in Fig.\ref{fig:CRNET} the total transmitted power, $\mathtt{TxP}$, is given by
\begin{align}
    \mathtt{TxP} & \triangleq \Exptn{\|\Mat{x}\|^2}  \nonumber \\
%                & = \Exptn{\Mat{x}^\dag \Mat{x}} = \Exptn{\Mat{s}^\dag\Mat{W}^\dag \Mat{W s}} \nonumber \\
%                & = \Tr{\Mat{W}^\dag \Mat{W}} = \Tr{\Mat{WW}^\dag} \nonumber \\
%                & = \sum_{k=1}^{K} \Tr{\Mat{w}_k \Mat{w}_k^\dag} = \sum_{k=1}^{K} \Tr{\Mat{w}_k^\dag \Mat{w}_k} \nonumber \\
                 & = \sum_{k=1}^{K} \| \Mat{w}_k \|^2 \label{eqn:TxP}.
\end{align}

The received signal at the $k$th SU is
\begin{align}
    y_k & = \Mat{h}_k^\dag \Mat{x} + n_k \nonumber \\
        & = \Mat{h}_k^\dag \Mat{w}_k s_k + \sum_{\substack{i=1\\i \neq k}}^{K} \Mat{h}_k^\dag \Mat{w}_i s_i + n_k. \label{eqn:SU_Rx_Power}
\end{align}
The right-hand side of \eqref{eqn:SU_Rx_Power} has three terms.
The first term is the received signal of the intended message, while the second and the third terms show the interference from other messages and noise, which is white and Gaussian, i.e. $n_k \sim \Set{CN}(0,\sigma_n^2)$, respectively.
The average received power for $k$th SU, $\mathtt{S}_k$, from the intended message is
\begin{align}
    \mathtt{S}_k & \triangleq \Exptn{| \Mat{h}_k^\dag \Mat{w}_k s_k |^2} \nonumber \\
%                & = \Exptn{(\Mat{h}_k \Mat{w}_k s_k)^\dag (\Mat{h}_k \Mat{w}_k s_k)} \nonumber \\
%                & = \Exptn{(\Mat{h}_k \Mat{w}_k s_k) (\Mat{h}_k \Mat{w}_k s_k)^\dag} \nonumber \\
%                & = \Mat{w}_k^\dag \Mat{h}_k \Mat{h}_k^\dag \Mat{w}_k \nonumber \\
                 & = | \Mat{w}_k^\dag \Mat{h}_k |^2 \label{eqn:S_k}.
\end{align}
Similarly, it is easy to show that the received interference power, $\mathtt{IP}_k$, is given by
\begin{align}
    \mathtt{IP}_k \triangleq \Exptn{\left|\sum_{\substack{i=1\\i\neq k}}^{K} \Mat{h}_k^\dag \Mat{w}_i s_i \right|^2} = \sum_{\substack{i=1\\i\neq k}}^{K}   |\Mat{w}_i^\dag \Mat{h}_k |^2. \label{eqn:IP_k}
\end{align}
Using \eqref{eqn:S_k} and \eqref{eqn:IP_k}, the SINR of $k$th SU-Rx, $\mathtt{SINR}_k$, is given by
\begin{align}
    \mathtt{SINR}_k \triangleq \frac{\mathtt{S}_k}{\sigma_n^2 + \mathtt{IP}_k} = \frac{| \Mat{w}_k^\dag \Mat{h}_k |^2}{\sigma_n^2 +  \sum_{\substack{i=1\\i\neq k}}^{K} |\Mat{w}_i^\dag \Mat{h}_k |^2}. \label{eqn:SINR_k}
\end{align}

In robust design problems relating to SINR, expressions with the form of $ |\Mat{w}_k^\dag \Mat{h}_k|^2 $ are frequently used. We can write
\begin{align}
	|\Mat{w}_k^\dag \Mat{h}_k|^2 & = \Mat{w}_k^\dag(\Mat{\tilde{h}}_k+\Mat{a}_k)(\Mat{\tilde{h}}_k+\Mat{a}_k)^\dag \Mat{w}_k \nonumber\\
%								   & = (\Mat{w}_k^\dag \Mat{h}_k) (\Mat{w}_k^\dag \Mat{h}_k)^\dag \nonumber\\
%								   & = \Mat{w}_k^\dag \Mat{h}_k \Mat{h}_k^\dag \Mat{w}_k \nonumber\\
								   & = \Mat{w}_k^\dag (\Mat{\tilde{H}}_k + \Mat{\Delta}_k) \Mat{w}_k,\label{eqn:qForm}
\end{align}
where $\Mat{\tilde{H}}_k = \Mat{\tilde{h}}_k \Mat{\tilde{h}}_k^\dag$ is a constant matrix and $\Mat{\Delta}_k$ is given by
\begin{equation}
  \Mat{\Delta}_k = \Mat{\tilde{h}}_k \Mat{a}_k^\dag + \Mat{a}_k \Mat{\tilde{h}}_k^\dag + \Mat{a}_k \Mat{a}_k^\dag.
\end{equation}
Note that $\Mat{\Delta}_k$ is a norm bounded matrix, $\|\Mat{\Delta}_k\| \leq \epsilon_k$.
It is straightforward to find the following relation:
\begin{align}
    \epsilon_k \geq \|\Mat{\Delta}_k\|&= \|\Mat{\tilde{h}}_k \Mat{a}_k^\dag + \Mat{a}_k \Mat{\tilde{h}}_k^\dag + \Mat{a}_k \Mat{a}_k^\dag\|\nonumber\\
	&\leq \|\Mat{\tilde{h}}_k \Mat{a}_k^\dag\| + \|\Mat{a}_k \Mat{\tilde{h}}_k^\dag\| + \|\Mat{a}_k \Mat{a}_k^\dag\|\nonumber\\
	&\leq \|\Mat{\tilde{h}}_k\|\ \|\Mat{a}_k^\dag\| + \|\Mat{a}_k\|\ \|\Mat{\tilde{h}}_k^\dag\| + \|\Mat{a}_k\|^2\nonumber\\
	&= \delta_k^2 + 2 \delta_k \|\Mat{\tilde{h}}_k\|. \label{eqn:sig_del_rel}
\end{align}
Using \eqref{eqn:sig_del_rel} it is possible to choose $\epsilon_k = \delta_k^2 + 2 \delta_k \|\Mat{\tilde{h}}_k\|$.
We may use the identity of $\Mat{x}^\dag\Mat{Ax} = \Tr{\Mat{Axx}^\dag}$, to further simplify this expression, which gives
\begin{align}
	|\Mat{w}_k^\dag \Mat{h}_k|^2 = \Tr{\left(\Mat{\tilde{H}}_k + \Mat{\Delta}_k\right) \Mat{w}_k\Mat{w}_k^\dag} .
\end{align}
We also adopt the notation of $\Mat{W}_k = \Mat{w}_k\Mat{w}^\dag_k$ in our design formulation.
With this, we find that
\begin{align}
	|\Mat{w}_k^\dag \Mat{h}_k|^2 = \Tr{\left(\Mat{\tilde{H}}_k + \Mat{\Delta}_k\right) \Mat{W}_k} .
\end{align}
It is noted that, from now on, similar expressions will appropriately be used for the other terms of $\mathtt{SINR}_k$.

Also, the received signal at the $\ell$th PU is
\begin{align}
    z_\ell & = \Mat{g}_\ell^\dag \Mat{x} + \nu_\ell \nonumber\\
           & = \sum_{k=1}^{K} \Mat{g}_\ell \Mat{w}_k s_k + \nu_\ell ,
\end{align}
and the interference power, $\mathtt{IP}_\ell$, to this  PU-Rx would be
\begin{align}
    \mathtt{IP}_\ell \triangleq \sum_{k=1}^{K} |\Mat{w}_k^\dag \Mat{g}_\ell |^2 . \label{eqn:IP_ell}
\end{align}
Again using similar formulation as $| \Mat{w}_k^\dag \Mat{h}_k|^2$, we get
\begin{align}
	| \Mat{w}_k^\dag \Mat{g}_\ell|^2 & = \Mat{w}_k^\dag (\Mat{\tilde{G}}_\ell + \Mat{\Lambda}_\ell) \Mat{w}_k \nonumber\\
	& = \Tr{\left(\Mat{\tilde{G}}_\ell + \Mat{\Lambda}_\ell \right) \Mat{W}_k} , \label{eqn:IP_term}
\end{align}
where $\Mat{\tilde{G}}_\ell$ is a constant matrix, $\Mat{\tilde{G}}_\ell = \Mat{\tilde{g}}_\ell \Mat{\tilde{g}}_\ell^\dag$ and $\Mat{\Lambda}_\ell$ is the norm bounded uncertainty matrix, $\|\Mat{\Lambda}_\ell\| \leq \xi_\ell$.
Similarly we know that $\xi_\ell = \eta_\ell^2 + 2 \eta_\ell \|\Mat{\tilde{g}}_\ell\|$.

Our design objective is to minimize the transmitted power, $\mathtt{TxP}$, while guaranteeing that the SINR at SU-Rx for all the channel realizations is higher than the QoS-constrained threshold, $\{\mathtt{SINR}_k \geq \gamma_i\}_{k=1}^{K}$, and simultaneously IP at PU-Rx is less than the PR-net--imposed threshold, $\{\mathtt{IP}_\ell \leq \kappa_\ell\}_{\ell=1}^{L}$, respectively.
Mathematically, this problem can be described as
\begin{align}
    \Min{\{\Mat{W}_k\}_{k=1}^{K}} & \quad \mathtt{TxP} \label{eqn:main_gen_form} \\
    \STo & \quad \mathop{\mathtt{SINR}_k}_{\forall \Mat{h}_k \in \Set{H}_k} \geq \gamma_k, \quad k = 1, \cdots, K \nonumber \\
         & \quad \mathop{\mathtt{IP}_\ell}_{\forall \Mat{g}_\ell \in \Set{G}_\ell} \leq \kappa_\ell, \quad \ell = 1, \cdots, L. \nonumber
\end{align}
The above problem is a problem with an infinite number of constraints.
To deal with such a problem one well-known method is to find the minimum and maximum values of $\mathtt{SINR}_k$ and $\mathtt{IP}_\ell$, respectively, related to those realizations of the channels which are claimed as the ``worst ones''.
The worst channel realizations for SINR and IP would lead to the minimum and maximum value of SINR and IP, respectively.
In that case, the problem will guarantee that the smallest possible SINR and largest possible IP also satisfy the constraints.
Using this worst case design methodology, we could recast \eqref{eqn:main_gen_form} to a simpler problem set as follows:
\begin{align}
    \Min{\{\Mat{W}_k\}_{k=1}^K}  & \quad \mathtt{TxP}    \label{eqn:symb_wc_form} \\
    \STo & \quad \min_{\Mat{h}_k \in \Set{H}_k}\ \mathtt{SINR}_k \geq \gamma_k, \quad k = 1, \cdots, K \nonumber \\
         & \quad \max_{\Mat{g}_\ell \in \Set{G}_\ell}\ \mathtt{IP}_\ell \leq \kappa_\ell, \quad \ell = 1, \cdots, L. \nonumber
\end{align}
Adopting the previously mentioned notations, we rewrite it as
\begin{subequations}\label{eqn:main_wc_prob}
\begin{align}
    \Min{\substack{\{\Mat{W}_k\}_{k=1}^{K}}} & \quad \sum_{k=1}^{K} \Tr{\Mat{W}_k}  \nonumber\\
%	\min_{\substack{\{\Mat{w}_k\}_{k=1}^{K}\\ \{\Mat{h}_k \in \Set{H}_k\}_{k=1}^{K} \\ \{\Mat{g}_\ell \in \Set{G}_\ell \}_{\ell=1}^{L}}} \sum_{k=1}^{K} \|\Mat{w}_k\|^2  & \quad \mathtt{subject\ to} \nonumber\\	
\STo  & \quad \min_{\|\Mat{\Delta}_k\| \leq \epsilon_k}\ \frac{\Tr{\left(\Mat{\tilde{H}}_k + \Mat{\Delta}_k\right)\Mat{W}_k}}{\sigma_n^2 +  \sum_{\substack{i=1\\i\neq k}}^{K} \Tr{\left(\Mat{\tilde{H}}_k + \Mat{\Delta}_k\right)\Mat{W}_i}} \geq \gamma_k, \quad k = 1, \cdots, K \label{eqn:SINRk_prob1} \\
    & \quad \max_{\|\Mat{\Lambda}_\ell\| \leq \eta_\ell}\ \sum_{k=1}^{K} \Tr{\left(\Mat{\tilde{G}}_\ell + \Mat{\Lambda}_\ell\right)\Mat{W}_k} \leq \kappa_\ell, \quad \ell = 1, \cdots, L . \label{eqn:IPell_prob1}
\end{align}
\end{subequations}
In the next sections, we will solve the robust problem of \eqref{eqn:main_wc_prob} and will show that this problem can be recast as a series of simple optimization problems.

%\section{Method 1: Conventional Robust Approximation \label{method1}}
\section{Loosely Bounded Robust Solution \label{method1}}
In this section we will deal with the problem of \eqref{eqn:main_wc_prob}.
In \cite{bengtsson1999} and \cite{bengtsson2001}, it is suggested to minimize the SINR through minimizing the numerator while maximizing its denominator.
So \eqref{eqn:SINRk_prob1} is equivalent to
\begin{align}
	& \min_{\| \Mat{\Delta}_k \| \leq \epsilon_k} \Tr{(\Mat{\tilde{H}}_k + \Mat{\Delta}_k) \Mat{W}_k } -
	\gamma_k \sum_{\substack{i=1 \\ i \neq k}}^{K} \max_{\| \Mat{\Delta}_k \| \leq \epsilon_k} \Tr{(\Mat{\tilde{H}}_k + \Mat{\Delta}_k) \Mat{W}_i} \geq \gamma_k \sigma_n^2 . \label{eqn:loose_min_sinr}
\end{align}
As it is known, this method is a conservative way to find the minimum of the SINR.

\subsection{Minimization of SINR}
To minimize the numerator,
\begin{align}
	\min_{\|\Mat{\Delta}_k\| \leq \epsilon_k}\ \Tr{\left(\Mat{\tilde{H}}_k + \Mat{\Delta}_k\right)\Mat{W}_k} ,
\end{align}
we adopt a loose lower bound, proposed by \cite{bengtsson1999}, \cite{bengtsson2001}.
Using this lower bound, we have
\begin{align}
	\min_{\|\Mat{\Delta}_k\| \leq \epsilon_k}\ \Tr{\left(\Mat{\tilde{H}}_k + \Mat{\Delta}_k\right)\Mat{W}_k} = \Tr{\left(\Mat{\tilde{H}}_k - \epsilon_k \Mat{I}_N\right)\Mat{W}_k} ,
\end{align}
and to maximize the denominator, the following term should be maximized
\begin{align}
	\max_{\|\Mat{\Delta}_k\| \leq \epsilon_k}\ \Tr{\left(\Mat{\tilde{H}}_k + \Mat{\Delta}_k\right)\Mat{W}_i}.
\end{align}
Using a similar approximation, we have
\begin{align}
	\max_{\|\Mat{\Delta}_k\| \leq \epsilon_k}\ \Tr{\left(\Mat{\tilde{H}}_k + \Mat{\Delta}_k\right)\Mat{W}_i} = \Tr{\left(\Mat{\tilde{H}}_k + \epsilon_k \Mat{I}_N\right)\Mat{W}_i}.
\end{align}
Using these results, the problem of SINR minimization \eqref{eqn:SINRk_prob1} is recast as
\begin{align}
	& \Tr{\left(\Mat{\tilde{H}}_k - \epsilon_k \Mat{I}_N\right)\Mat{W}_k} - \gamma_k \sum_{\substack{i=1\\i \neq k}}^{K} \Tr{\left(\Mat{\tilde{H}}_k + \epsilon_k \Mat{I}_N\right)\Mat{W}_i} \geq \sigma_n^2\gamma_k, \ k = 1, \cdots, K ,
\end{align}
and by regrouping the left hand side of this equation we find
\begin{align}
	& \Tr{\Mat{\tilde{H}}_k \left(\Mat{W}_k - \gamma_k \sum_{\substack{i=1\\i\neq k}}^{K} \Mat{W}_i \right)} - \epsilon_k \Tr{\Mat{W}_k + \gamma_k \sum_{\substack{i=1\\i\neq k}}^{K} \Mat{W}_i } \geq \sigma_n^2\gamma_k, \ k = 1, \cdots, K.
\end{align}

\subsection{The Whole Conventional Program}
Using the same methodology as before, IP maximization \eqref{eqn:IPell_prob1} leads to the following problem
\begin{align}
	\sum_{k=1}^{K}\Tr{\left(\Mat{\tilde{G}}_\ell + \xi_\ell \Mat{I}_N \right)\Mat{W}_k} \leq \kappa_\ell, \quad \ell = 1, \cdots, L.
\end{align}

Then the whole conventional program targeting to solve the robust downlink optimization in MISO CR-Nets becomes
\begin{subequations}\label{eqn:method1}
\begin{align}
    \Min{\substack{\{\Mat{W}_k\}_{k=1}^{K}}} & \quad \sum_{k=1}^{K} \Tr{\Mat{W}_k} \nonumber \\
	\STo & \quad\Tr{\Mat{\tilde{H}}_k \left(\Mat{W}_k - \gamma_k \sum_{\substack{i=1\\i\neq k}}^{K} \Mat{W}_i \right)} - \epsilon_k \Tr{\Mat{W}_k + \gamma_k \sum_{\substack{i=1\\i\neq k}}^{K} \Mat{W}_i } \geq \sigma_n^2\gamma_k, \ k = 1, \cdots, K  \\
	& \quad \Tr{\left(\Mat{\tilde{G}}_\ell + \xi_\ell \Mat{I}_N \right) \sum_{k=1}^{K}\Mat{W}_k} \leq \kappa_\ell, \quad \ell = 1, \cdots, L, \\
	& \qquad \Mat{W}_k = \Mat{W}_k^\dag, \quad k=1, \cdots, K,\\
	& \qquad \Mat{W}_k \succeq 0, \quad k=1, \cdots, K.
\end{align}
\end{subequations}
Please note the fact that the last two constraints are inherent in the structure of the problem formulation.
Also note that to come up with a convex problem formulation, a non-convex constraint, $\mathsf{rank}\{\Mat{W}_k\} = 1$, is eliminated \cite{vorobyov2003}, \cite{bengtsson1999}, \cite{sharma2008}.
This final form of the problem is an SDP and can be solved using efficient numerical methods \cite{cvx2006}.
Finally it should be noted that unlike \cite{{sharma2008}} the beamforming weights are not exactly the principal eigenvector\footnote{The principal eigenvector of a rank one matrix is the eigenvector corresponding to the only non-zero eigenvalue.} of the matrix solution.
To get the beamforming weights, the eigen decomposition of the $\Mat{W}_k$ is used.
In this decomposition, $\Mat{W}_k$ may be decomposed to a series of rank one matrices, i.e.,
\begin{align}
	\Mat{W}_k = \sum_{n=1}^{N} \lambda_{n,k}\ \Mat{e}_{n,k}\ \Mat{e}_{n,k}^\dag,	
\end{align}
where in this expansion, $\lambda_{n,k}$ denotes the $n$th eigenvalue and $\Mat{e}_{n,k}$ is its respective eigenvector.
The solution matrix of $\Mat{W}_k$ itself is a rank one matrix, then all the eigenvalues are equal to zero except one, let's say $\lambda_{N,k}$.
Therefore the above mentioned equation may be written as
\begin{align}
	\Mat{W}_k & = \lambda_{N,k}\ \Mat{e}_{N,k}\ \Mat{e}_{N,k}^\dag	\nonumber \\
			  & = (\sqrt{\lambda_{N,k}}\ \Mat{e}_{N,k}) (\sqrt{\lambda_{N,k}}\ \Mat{e}_{N,k})^\dag \nonumber \\
			  & = \Mat{w}_k \Mat{w}_k^\dag,
\end{align}
where $\Mat{w}_k = \sqrt{\lambda_{N,k}}\ \Mat{e}_{N,k}$.
The next two sections deal with the same problem but in different ways.
It should be noted that the next two problems are not SDP and are generally convex problems, but, for these problems we use the same formula to acquire the beamforming weights from the solution matrix.

%\section{Method 2: Less-Conservative Robust Approximation \label{method2}}
\section{Strictly Bounded Robust Solution \label{method2}}
In the previous section, the minimum of SINR was found using loose upper and lower bounds for its constituent terms.
In this section, we try to minimize the SINR using the same method: we minimize the numerator and maximizing the denominator.
But here, we try to find the exact maximum and the exact minimum for each term respectively:
\begin{align}
	& \min_{\| \Mat{\Delta}_k \| \leq \epsilon_k} \Tr{(\Mat{\tilde{H}}_k + \Mat{\Delta}_k) \Mat{W}_k } - \gamma_k \sum_{\substack{i=1 \\ i \neq k}}^{K} \max_{\| \Mat{\Delta}_k \| \leq \epsilon_k} \Tr{(\Mat{\tilde{H}}_k + \Mat{\Delta}_k) \Mat{W}_i} \geq \gamma_k \sigma_n^2 . \label{eqn:strict_min_sinr}
\end{align}
Our main tool, is the Lagrangian Multiplier method.

\subsection{Minimization of SINR}
We start with the first minimization problem.

{\em Proposition 1:} For the terms $\Tr{\left(\Mat{\tilde{H}}_k + \Mat{\Delta}_k\right) \Mat{W}_k}$, using a norm-bounded variable $\Mat{\Delta}_k$, $\|\Mat{\Delta}_k\| \leq \epsilon_k$, the minimizer and maximizer would be
\begin{equation}
	\Mat{\Delta}_k^{min} = - \epsilon_k \frac{\Mat{W}_k^\dag}{\|\Mat{W}_k\|},
\end{equation}
and
\begin{equation}
	\Mat{\Delta}_k^{max} = \epsilon_k \frac{\Mat{W}_k^\dag}{\|\Mat{W}_k\|},
\end{equation}
respectively.

{\em Proof:} Please refer to Appendix \ref{app1}.

Using the above results, we have
\begin{align}
	& \min_{\| \Mat{\Delta}_k \| \leq \epsilon_k} \Tr{(\Mat{\tilde{H}}_k + \Mat{\Delta}_k) \Mat{W}_k } = \Tr{\left(\Mat{\tilde{H}}_k - \epsilon_k \frac{\Mat{W}_k^\dag }{\| \Mat{W}_k \|} \right) \Mat{W}_k} =  \Tr{\Mat{\tilde{H}}_k \Mat{W}_k} - \epsilon_k \|\Mat{W}_k\| , \\
	& \max_{\| \Mat{\Delta}_k \| \leq \epsilon_k} \Tr{(\Mat{\tilde{H}}_k + \Mat{\Delta}_k) \Mat{W}_i } = \Tr{\left(\Mat{\tilde{H}}_k + \epsilon_k \frac{\Mat{W}_i^\dag}{\| \Mat{W}_i \|} \right) \Mat{W}_i} =  \Tr{\Mat{\tilde{H}}_k \Mat{W}_i} + \epsilon_k \|\Mat{W}_i\| .
\end{align}
So we may rewrite (\ref{eqn:strict_min_sinr}) as
\begin{align}
 & \Tr{\Mat{\tilde{H}}_k \Mat{W}_k} - \epsilon_k \|\Mat{W}_k \| - \gamma_k \sum_{\substack{i=1\\i\neq k}}^{k}  \left( \Tr{\Mat{\tilde{H}}_k \Mat{W}_i} + \epsilon_k \|\Mat{W}_i \| \right) \\
 & \ = \Tr{\Mat{\tilde{H}}_k \left( \Mat{W}_k - \gamma_k \sum_{\substack{i=1\\i\neq k}} \Mat{W}_i \right) } - \epsilon_k \left(\|\Mat{W}_k\| + \gamma_k \sum_{\substack{i=1\\i\neq k}} \|\Mat{W}_i\|\right) \geq \gamma_k \sigma_n^2 .
\end{align}

\subsection{The Whole Program}
Similarly, the IP constraints may be written as:
\begin{align}
	\max_{\| \Mat{\Lambda}_\ell \| \leq \xi_k} & \sum_{k=1}^{K} \Tr{(\Mat{\tilde{G}}_\ell + \Lambda_\ell) \Mat{W}_k} = \sum_{k=1}^{K} \left( \Tr{\Mat{\tilde{G}}_\ell \Mat{W}_k} + \xi_\ell \|\Mat{W}_k \| \right) \leq \kappa_\ell .
\end{align}
Finally, the whole program is
\begin{subequations}\label{eqn:method2}
\begin{align}
    \Min{\substack{\{\Mat{W}_k\}_{k=1}^{K}}}& \quad \sum_{k=1}^{K} \Tr{\Mat{W}_k} \nonumber \\
	\STo & \quad\Tr{\Mat{\tilde{H}}_k \left(\Mat{W}_k - \gamma_k \sum_{\substack{i=1\\i\neq k}}^{K} \Mat{W}_i \right)} - \epsilon_k \left( \|\Mat{W}_k\| + \gamma_k \sum_{\substack{i=1\\i\neq k}}^{K} \|\Mat{W}_i\| \right) \geq \sigma_n^2\gamma_k, \nonumber \\
	& \qquad \qquad \qquad k = 1, \cdots, K  \\
	& \quad \sum_{k=1}^{K} \left( \Tr{\Mat{\tilde{G}}_\ell \Mat{W}_k} + \xi_\ell \|\Mat{W}_k\| \right) \leq \kappa_\ell, \quad \ell = 1, \cdots, L \\
	& \quad \Mat{W}_k = \Mat{W}_k^\dag, \quad k=1, \cdots, K \\
	& \quad \Mat{W}_k \succeq 0, \quad k=1, \cdots, K .
\end{align}
\end{subequations}
Although this final problem is not an SDP, it is in fact convex, and this problem can be solved using standard numerical optimization packages, like CVX \cite{cvx2006}.

\section{Exact Robust Solution \label{method3}}
In the last two sections we dealt with the problem of minimizing the SINR using a conservative method.
In this section we find the exact solution instead.
We start again with the problem of \eqref{eqn:main_wc_prob}, but with a simple alteration.
This problem is stated as:
\begin{subequations}\label{eqn:General_Main_Problem}
	\begin{align}
		\Min{\{\Mat{W}_k\}_{k=1}^{K}} & \quad \sum_{k=1}^{K} \Tr{\Mat{W}_k} \nonumber \\
		\STo & \ \min_{\|\Mat{\Delta}_k\| \leq \epsilon_k} \left( \Tr{(\Mat{\tilde{H}}_k + \Mat{\Delta}_k) \Mat{W}_k} -	\gamma_k \sum_{\substack{i=1\\i\neq k}}^{K} \Tr{(\Mat{\tilde{H}}_k+\Mat{\Delta}_k) \Mat{W}_i} \right) \geq \sigma_n^2 \gamma_k, \ \ k=1, \cdots, K ; \label{eqn:exact_min_sinr}\\
		& \ \max_{\|\Mat{\Lambda}_\ell\| \leq \xi_\ell} \sum_{k = 1}^{K} \Tr{(\Mat{\tilde{G}}_\ell + \Mat{\Lambda}_\ell) \Mat{W}_k} \leq \kappa_\ell, \quad \ell = 1, \cdots, L .
	\end{align}
\end{subequations}
In the above problem, we try to minimize the SINR directly and without using conservative assumptions.
First we have the following proposition:

{\em Proposition 2:} The minimizer of \eqref{eqn:exact_min_sinr} has the form of
\begin{equation}
	\Mat{\Delta}_k^{min} = -\epsilon_k \frac{\left(\Mat{W}_k - \gamma_k \sum_{\substack{i=1\\ i\neq k}}^{K} \Mat{W}_i\right)^\dag}{\|\Mat{W}_k - \gamma_k \sum_{\substack{i=1\\ i\neq k}}^{K} \Mat{W}_i\|} .
\end{equation}

{\em Proof:} Please refer to Appendix \ref{app2}.

Using this proposition we find
\begin{align}
	\min_{\|\Mat{\Delta}_k\| \leq \epsilon_k} \left( \Tr{(\Mat{\tilde{H}}_k + \Mat{\Delta}_k) \Mat{W}_k} - \gamma_k \sum_{\substack{i=1\\i\neq k}}^{K} \Tr{(\Mat{\tilde{H}}_k+\Mat{\Delta}_k) \Mat{W}_i} \right) & = \nonumber \\
	% & = \Tr{(\Mat{\tilde{H}}_k + \Mat{\Delta}_k^*) \Mat{W}_k} - \nonumber \\
	% & \qquad \qquad \gamma_k \sum_{\substack{i=1\\i\neq k}}^{K} \Tr{(\Mat{\tilde{H}}_k+\Mat{\Delta}_k^*) \Mat{W}_i}  \nonumber \\
	% & = \Tr{\Mat{\tilde{H}}_k \Mat{W}_k} - \gamma_k \sum_{\substack{i=1\\i\neq k}}^{K} \Tr{\Mat{\tilde{H}}_k \Mat{W}_i} + \nonumber \\
	% & \qquad \qquad\Tr{\Mat{\Delta}_k^* \Mat{W}_k} - \gamma_k \sum_{\substack{i=1\\i\neq k}}^{K} \Tr{\Mat{\Delta}_k^* \Mat{W}_i}  \nonumber \\
	% & = \Tr{\Mat{\tilde{H}}_k \Mat{W}_k} - \gamma_k \sum_{\substack{i=1\\i\neq k}}^{K} \Tr{\Mat{\tilde{H}}_k \Mat{W}_i} + \nonumber \\
	% & \qquad \qquad\Tr{\Mat{\Delta}_k^* \left( \Mat{W}_k - \gamma_k \sum_{\substack{i=1\\i\neq k}}^{K} \Mat{W}_i \right) }  \nonumber \\
	\Tr{\Mat{\tilde{H}}_k ( \Mat{W}_k - \gamma_k \sum_{\substack{i=1\\i\neq k}}^{K} \Mat{W}_i ) } - \epsilon_k\ \| \Mat{W}_k - \gamma_k \sum_{\substack{i=1\\i\neq k}}^{K} \Mat{W}_i \| & \geq \sigma_n^2\gamma_k .
\end{align}

Finally we come up with the final and general problem:
\begin{subequations}\label{eqn:method3}
\begin{align}
	\Min{\{\Mat{W}_k\}_{k=1}^{K}} & \quad \sum_{k=1}^{K} \Tr{\Mat{W}_k} \nonumber \\
	\STo & \qquad \Tr{\Mat{\tilde{H}}_k ( \Mat{W}_k - \gamma_k \sum_{\substack{i=1\\i\neq k}}^{K} \Mat{W}_i ) } - \epsilon_k\ \| \Mat{W}_k - \gamma_k \sum_{\substack{i=1\\i\neq k}}^{K} \Mat{W}_i \| \geq \sigma_n^2\gamma_k, \quad k=1, \cdots, K;  \\
		& \qquad \sum_{k=1}^{K} \left( \Tr{\Mat{\tilde{G}}_\ell \Mat{W}_k} + \xi_\ell \|\Mat{W}_k \| \right) \leq \kappa_\ell, \quad \ell = 1, \cdots, L;\\
		& \qquad \Mat{W}_k = \Mat{W}_k^\dag, \quad k = 1, \cdots, K  \\
		& \qquad \Mat{W}_k \succeq 0,  \quad k = 1, \cdots, K .
	\end{align}
\end{subequations}
The above problem is the most general form of the original problem.
It should be noted that the beamforming weights are also the principal eigenvector of the solutions of this problem.
Also it should be mentioned that the IP part of these two last problems, \eqref{eqn:method2} and \eqref{eqn:method3}, are the same, i.e., these two problems have the same performance in IPs.

\section{Simulation Results and Discussions\label{simResults}}
To validate our developed methods, a set of simulations were conducted.
It is assumed that the BS is equipped with a Uniform Linear Array (ULA) having 8 elements with a spacing of half wave length.
A set of $K=3$ SU-Rx's are served and the CR-Net should protect a set of $L=2$ PU-Rx's.
The SU-Rx's are located in the directions of $\theta_1=20^\circ$, $\theta_2=35^\circ$ and $\theta_3=50^\circ$ relative to the antenna boreside, respectively.
The PU-Rx's are also located at the directions of $\phi_1 = 80^\circ$ and $\phi_2 = 85^\circ$, respectively.
It is assumed that the change in Direction of Arrival (DoA) of input waves to the SU-Tx may be changed up to $\pm 5^\circ$ arbitrarily.
The noise power is assumed to be $\sigma_n^2 = 0.01$, and constant for all of the users.
Also, a constant SINR level of 10dB is targeted for all the SUs, while the constant interference threshold of 0.01 is used to protect the PUs.
The channel model for PUs and SUs is assumed to be in line with a simple model of
\begin{align}
	[\Mat{h}_k(\theta_k)]_{i} &= e^{j \pi (i-1) \cos(\theta_k)}, \quad i,k=1,\cdots,K, \\
	[\Mat{g}_\ell(\phi_\ell)]_{i} &= e^{j \pi (i-1) \cos(\phi_\ell)}, \quad i,\ell=1,\cdots,L.
\end{align}
The uncertainty sets are characterized with $\epsilon_k = \eta_\ell = 0.05$.
We have used the CVX Software Package \cite{cvx2006} to solve the proposed problems numerically.

\begin{figure}[p]
	\centering
	\begin{tabular}{cc}
		\includegraphics[width=0.50\textwidth]{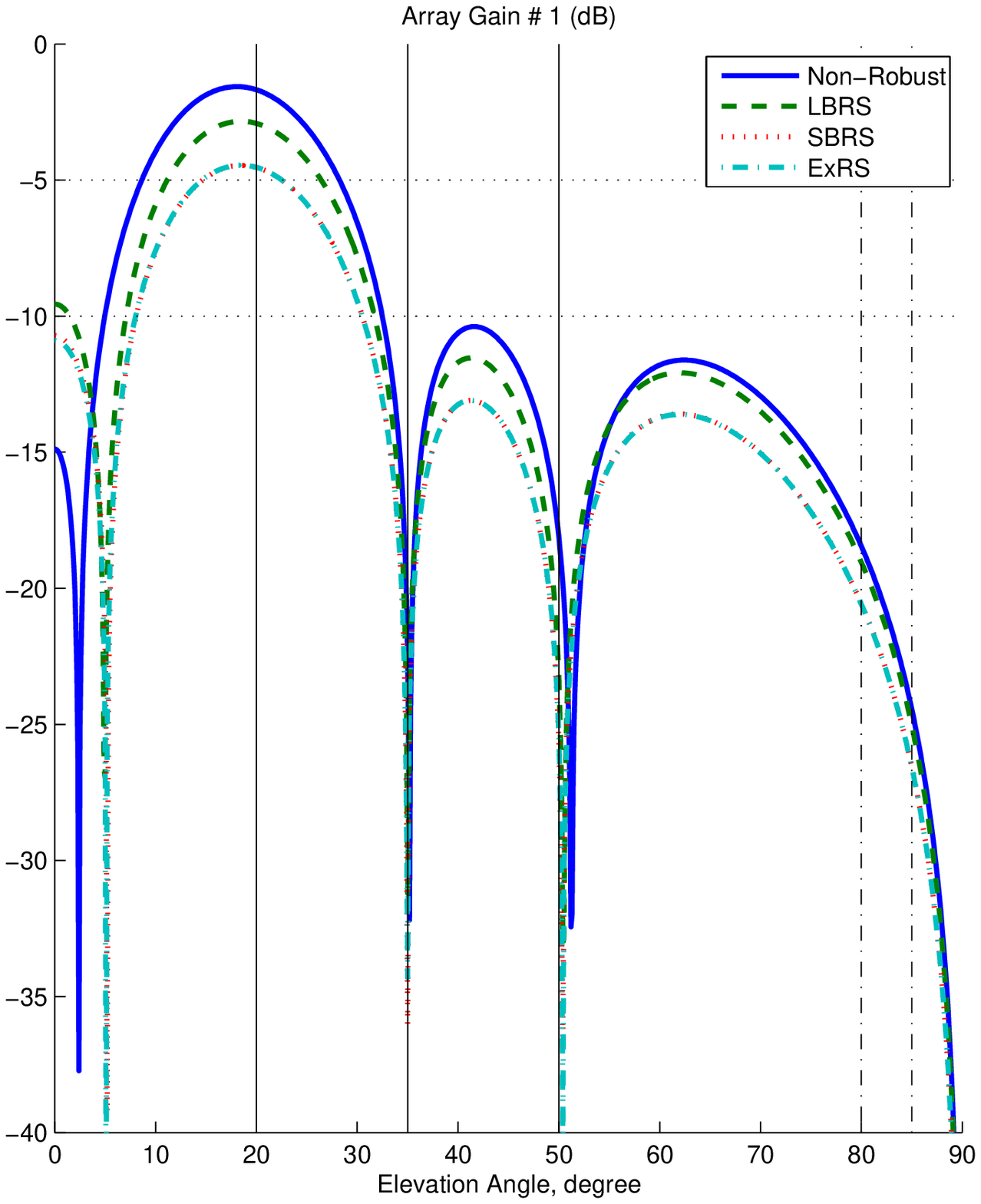} &
		\includegraphics[width=0.50\textwidth]{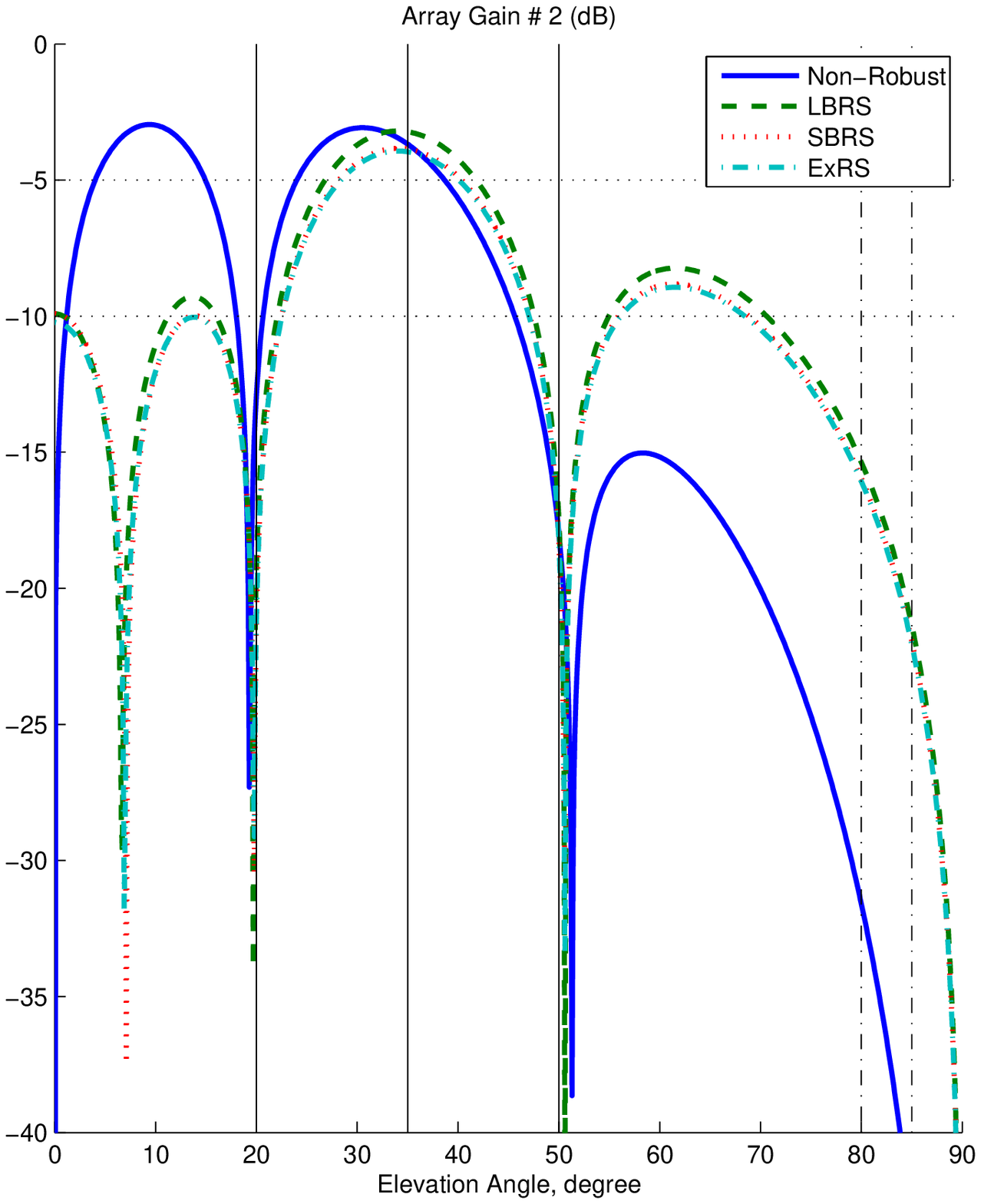} \vspace{-0.8cm}\\
		\scriptsize (a) Array Gain for SU\#1	& 	
		\scriptsize (b) Array Gain for SU\#2
	\end{tabular}
	\begin{tabular}{c}
		\includegraphics[width=0.50\textwidth]{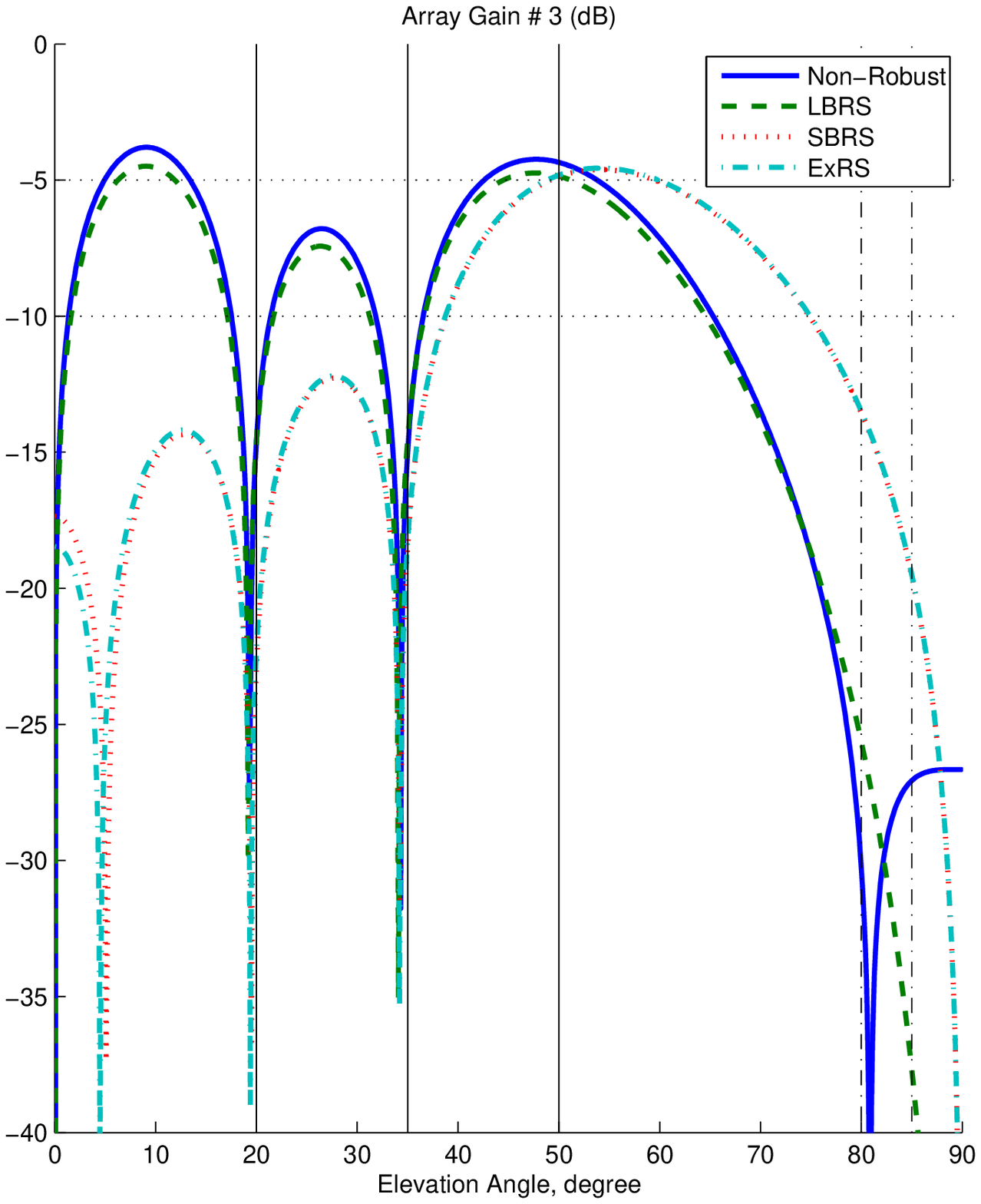} \vspace{-0.8cm}\\
		{\scriptsize (c) Array Gain for SU\#3}
	\end{tabular}
	\caption{Array Gain for Different Users}
	\label{fig:arrayGain}
\end{figure}

In Fig.~\ref{fig:arrayGain}, the array gains toward each user using different weight vectors are depicted.
In this figure and the subsequent figures, LCBS, SCBS and ExCS denote Loosely Bounded Convex Solution, Strictly Bounded Convex Solution, and Exact Convex Solution, respectively.
In this figure, the vertical solid lines show the DoA corresponding to different SU-Rx's, while the vertical dashed lines show the DoA of PU-Rx's.
From this figure, it is clear that the proposed method can transmit the desired data to the SU-Rx's while it is protecting the PU-Rx's. It is also apparent that all these approaches produce similar results.

\begin{figure}[p]
	\centering
	\begin{tabular}{cc}
		\includegraphics[width=0.45\textwidth]{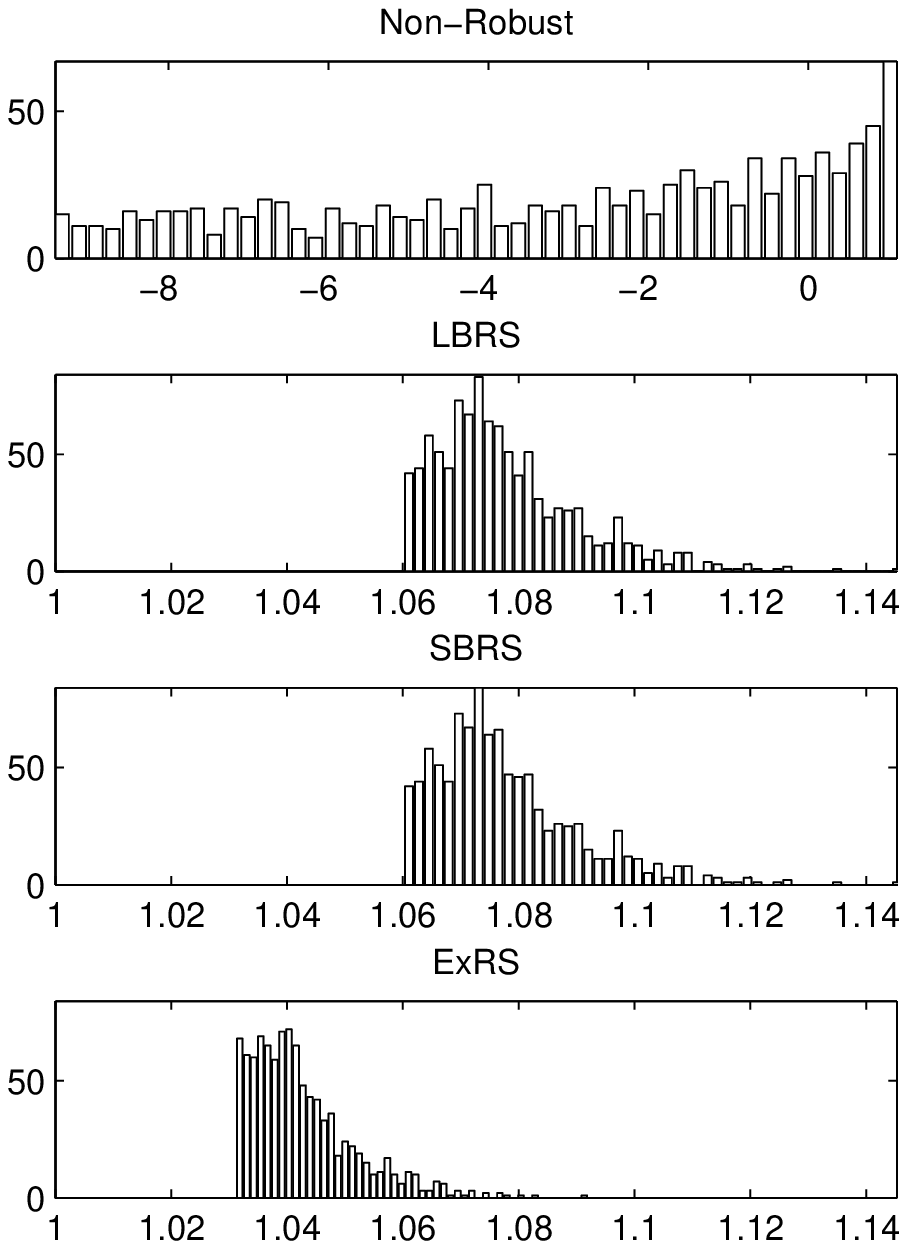} 	&
		\includegraphics[width=0.45\textwidth]{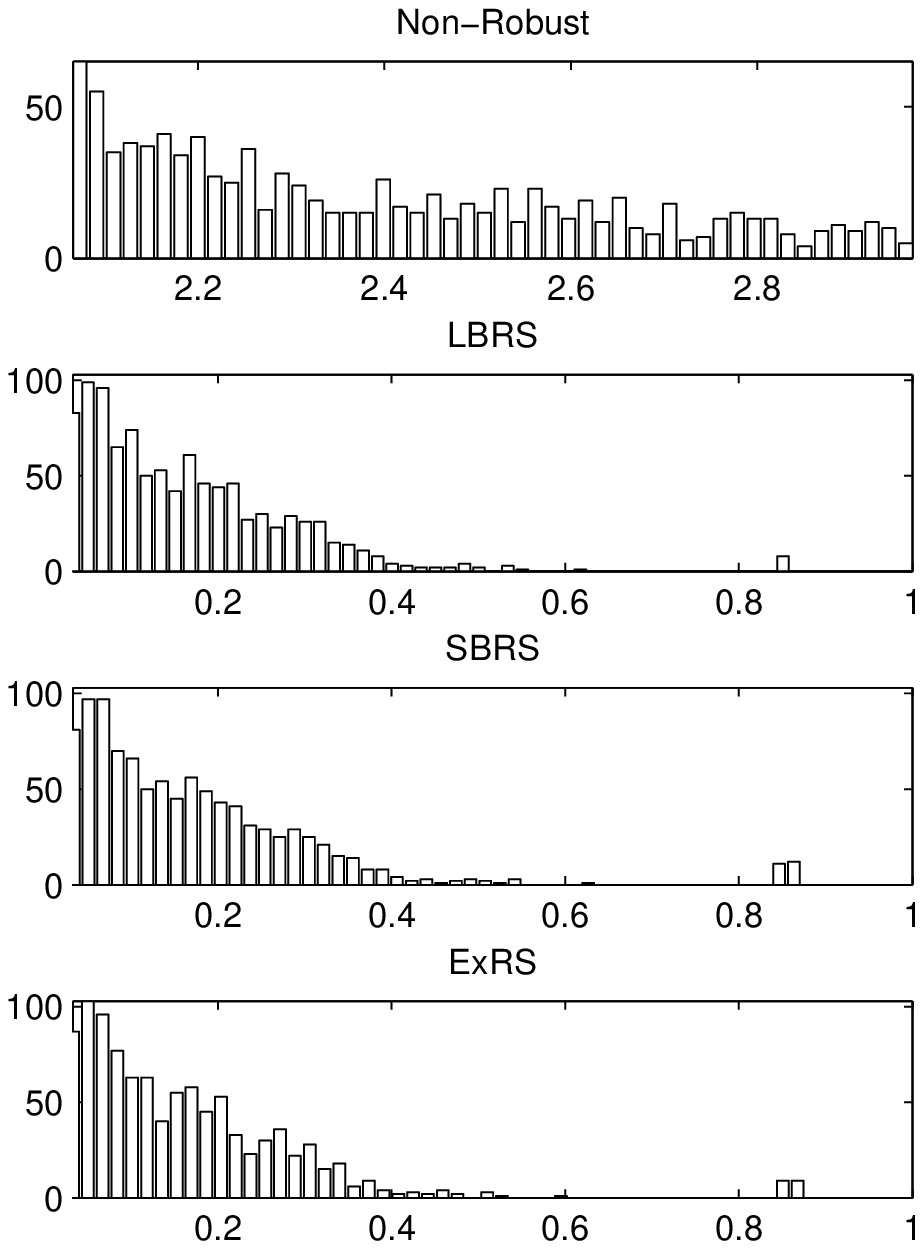}   	\vspace{-1cm}\\
		\scriptsize (a) Normalized SINR Constraints for SU\#1, $\epsilon_k = \xi_\ell = 0.05$		&
		\scriptsize (b) Normalized IP Constraints for PU\#1, $\epsilon_k = \xi_\ell = 0.05$			\\
		\includegraphics[width=0.45\textwidth]{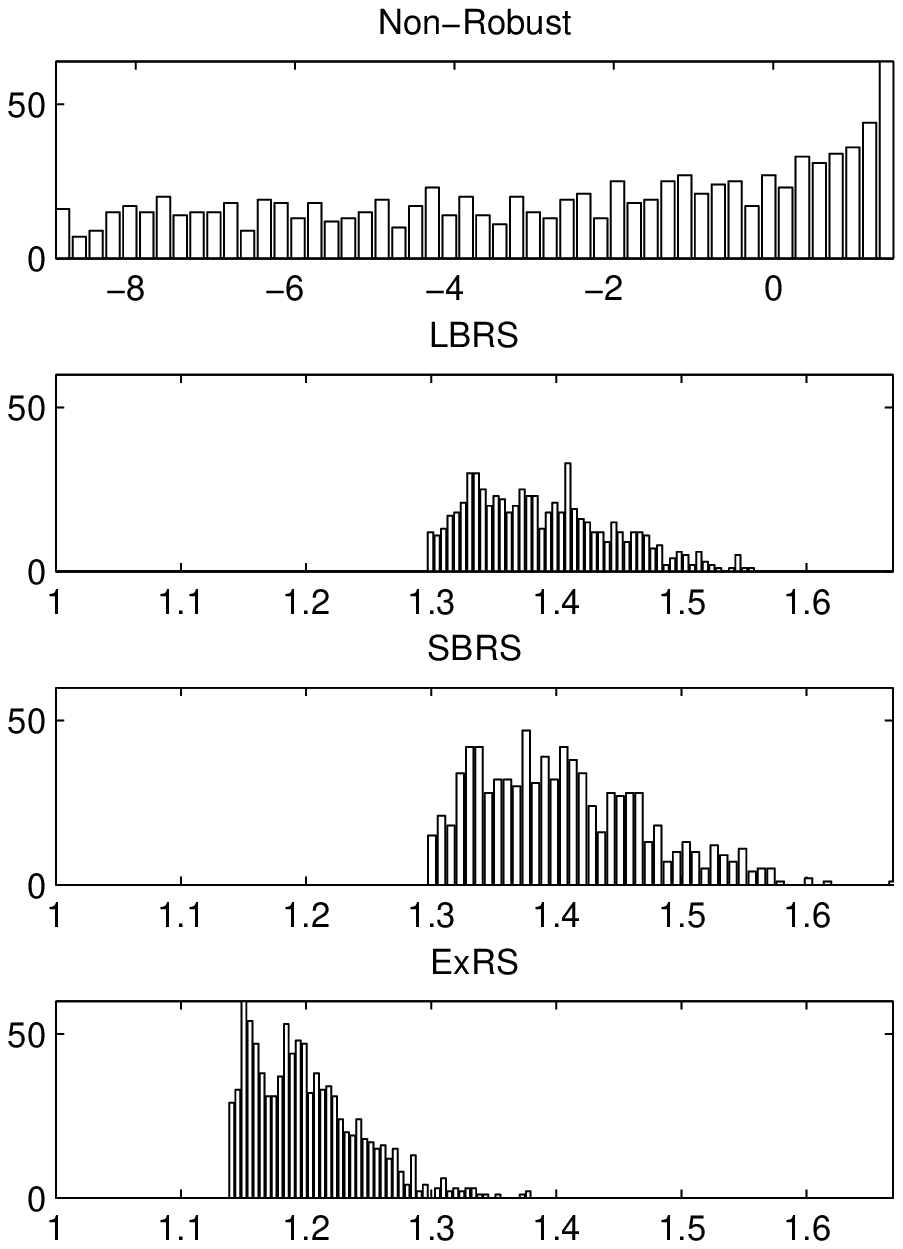}	&
		\includegraphics[width=0.45\textwidth]{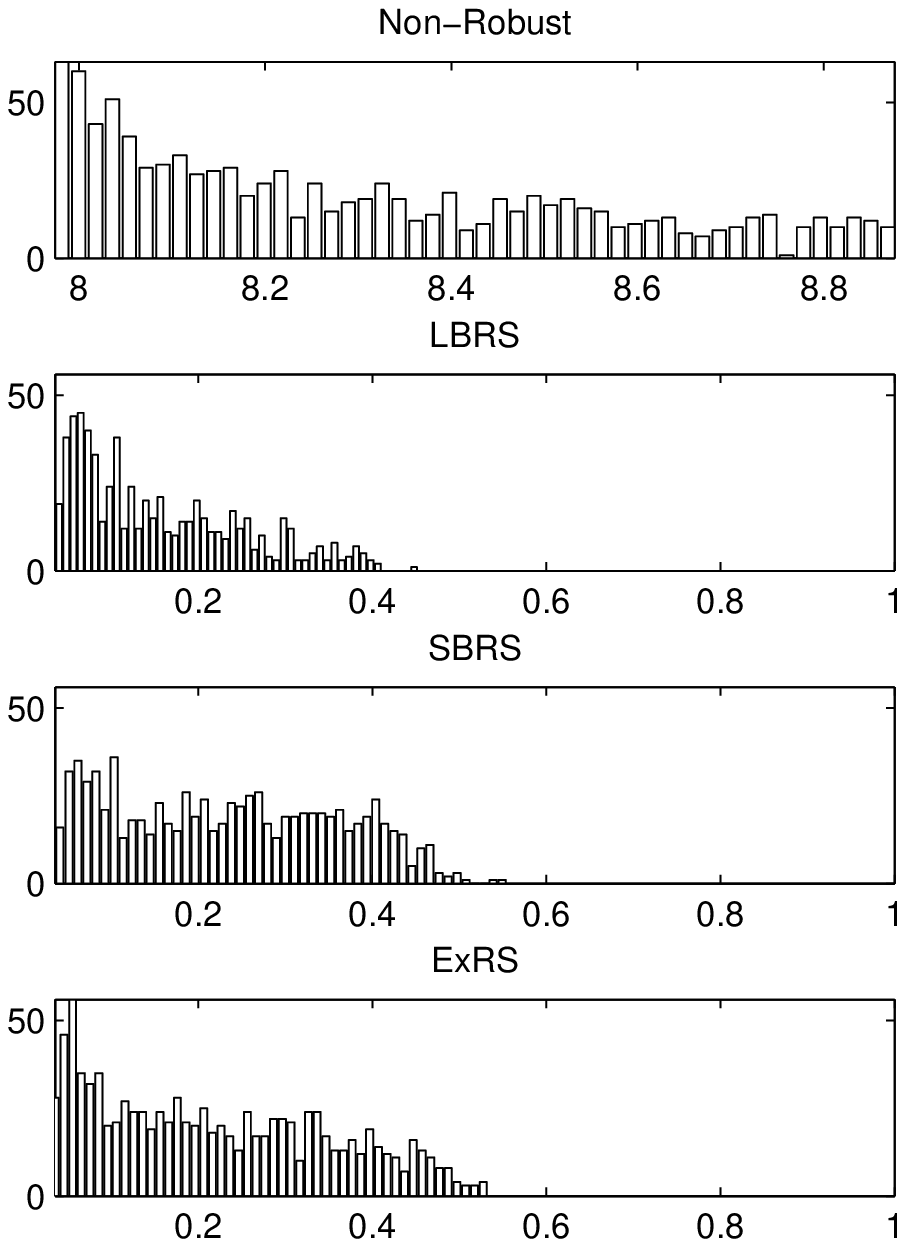}	\vspace{-1cm}\\
		\scriptsize (c) Normalized SINR Constraints for SU\#1, $\epsilon_k = \xi_\ell = 0.20$		&
		\scriptsize (d) Normalized IP Constraints for PU\#1, $\epsilon_k = \xi_\ell = 0.20$			
	\end{tabular}
	\caption{Normalized SINR Constraints for Different Methods for SU\#1 and PU\#1 }
	\label{fig:Hist_SU3_all}
\end{figure}

In Fig.~\ref{fig:Hist_SU3_all}, we have plotted the histogram of normalized constraints for both SU-Rx's and PU-Rx's.
The normalized constraints of SUs, $C_{k}^{\mathrm{(sinr)}}$, defined as ${\mathtt{SINR}_k}/{\gamma_k}$, is equivalent to
\begin{align}
	C_{k}^{\mathrm{(sinr)}} = \frac{1}{\sigma_n^2 \gamma_k}\Mat{w}_k^\dag \Mat{H}_k \Mat{w}_k - \frac{1}{\sigma_n^2} \sum_{\substack{i=1\\i\neq k}}^{K} \Mat{w}_i^\dag \Mat{H}_k \Mat{w}_i \nonumber ,
\end{align}
where $\Mat{H}_k = \Mat{h}_k\Mat{h}_k^\dag$, and the normalized PU constraint, $C_{k}^{\mathrm{(ip)}}$, is defined as
\begin{align}
	C_{\ell}^{\mathrm{(ip)}} = \frac{1}{\kappa_\ell} \sum_{k=1}^{K} \Mat{g}^\dag_\ell \Mat{H}_k \Mat{g}_\ell . \nonumber
\end{align}
Unlike the normalized SINR constraints for a SU-Rx, when a normalized IP constraint is less than one, this constraint is considered to be satisfied.
In Fig.~\ref{fig:Hist_SU3_all}-a and Fig.~\ref{fig:Hist_SU3_all}-c, the normalized SINR histograms for two different scenarios are depicted.
In the first one, the uncertainty sets are chosen to be $\epsilon_k = \xi_\ell = 0.05$ while for the second scenario, the uncertainty measures are four times more than the first one, i.e. $\epsilon_k = \xi_\ell = 0.20$.
It is apparent that in this scenario the gap between the ultimate value of normalized SINR, i.e. 1, and the actual values is wider than in the first scenario, having smaller uncertainty sets.
In both cases, as expected, ExCS is outperforming the other two schemes due to its exact bounds on SINR.
In Fig.~\ref{fig:Hist_SU3_all}-b and Fig.~\ref{fig:Hist_SU3_all}-d the normalized IP constraints for PU-Rx's are depicted.
As can be seen, there is little difference between the proposed methods in terms of their IP constraints, because of the similar structure of these constraints.
%\begin{figure}
%	\centering
%	\begin{tabular}{cc}
%	\includegraphics[clip,width=0.50\textwidth]{./img/Hist_SU1_all.eps} & \hspace{-1.2cm}
%	\includegraphics[clip,width=0.55\textwidth]{./img/Hist_SU1_all_x4.eps} \vspace{-1cm}\\
%		test \#1 & test \#2 \\
%	\includegraphics[width=0.45\textwidth]{./img/Hist_PU1_all.eps} \\
%	test \#3
%%	\epsfig{file=figure4.eps,width=0.5\linewidth,clip=}
%	\end{tabular}
%	\caption{test Multi fig}
%\end{figure}
As it can be seen, the robust design is immune for the variation of channel, whereas the non-robust model fails in such situations.
Also it is apparent that for the robust case, the variation of normalized constraints is much less than the variation of normalized constraints for the non-robust case.
It is also clear that SBCS and ExCS are more efficient in terms of handling the SINR.
The ExCS model not only can satisfy all the constraints, but also is a parsimonious model in terms of SINR.
This is because of the fact that this method uses exact minimum of SINR.
Also it should be noted that the IP variations in SBCS and ExCS are the same.
It is because of the identical form of the equation which describes the IP constraints in these two methods.
Additionally, it should be mentioned that they are slightly better than the LCBS.

\begin{figure}[p]
	\centering
	\begin{tabular}{c}
		\includegraphics[width=0.63\textwidth]{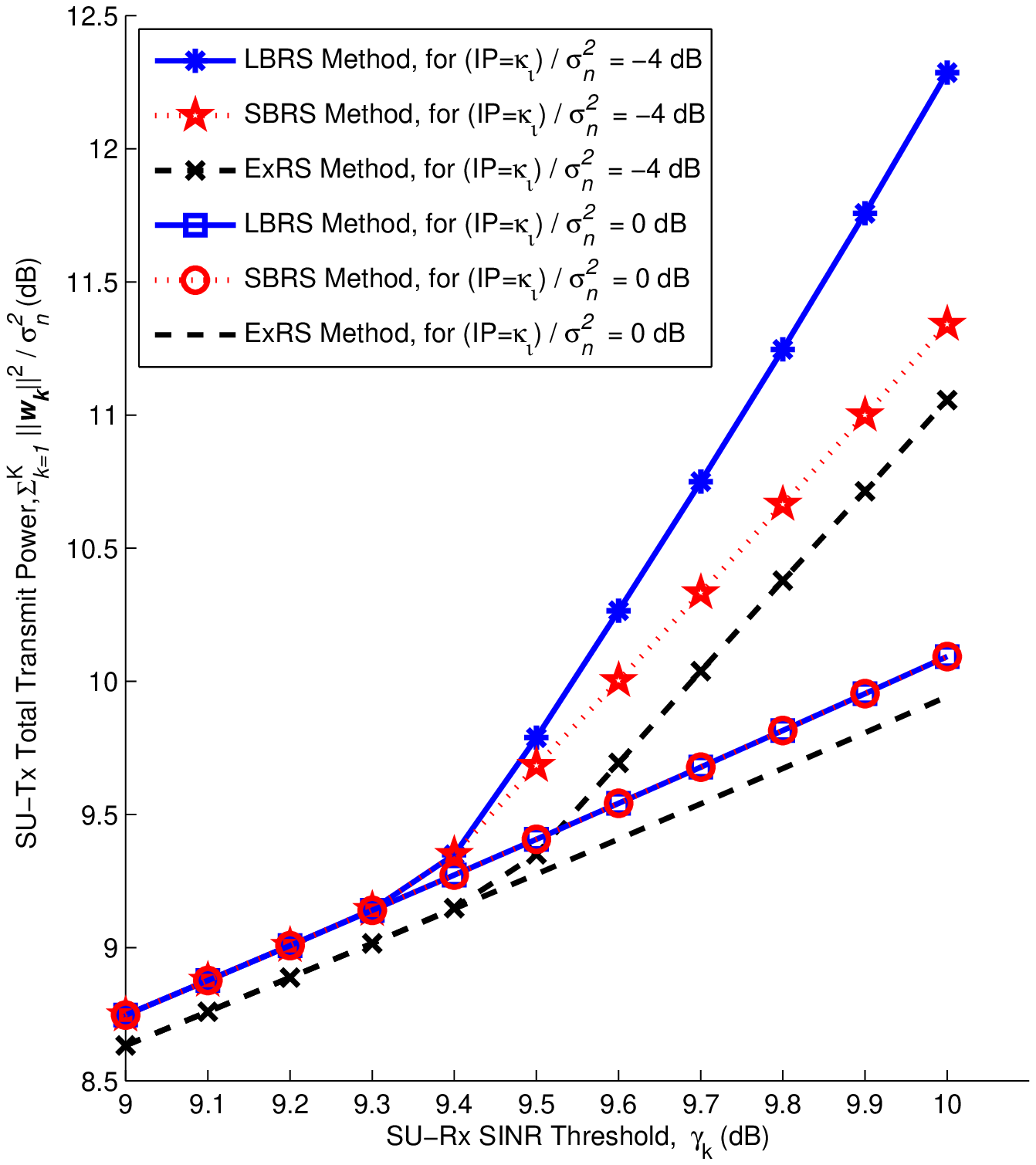} 	\vspace{-0.8cm}\\
		\scriptsize (a) Relative IP Levels of $-$4 dB and 0 dB, in Low SNR Regime		\\					
		\includegraphics[width=0.63\textwidth]{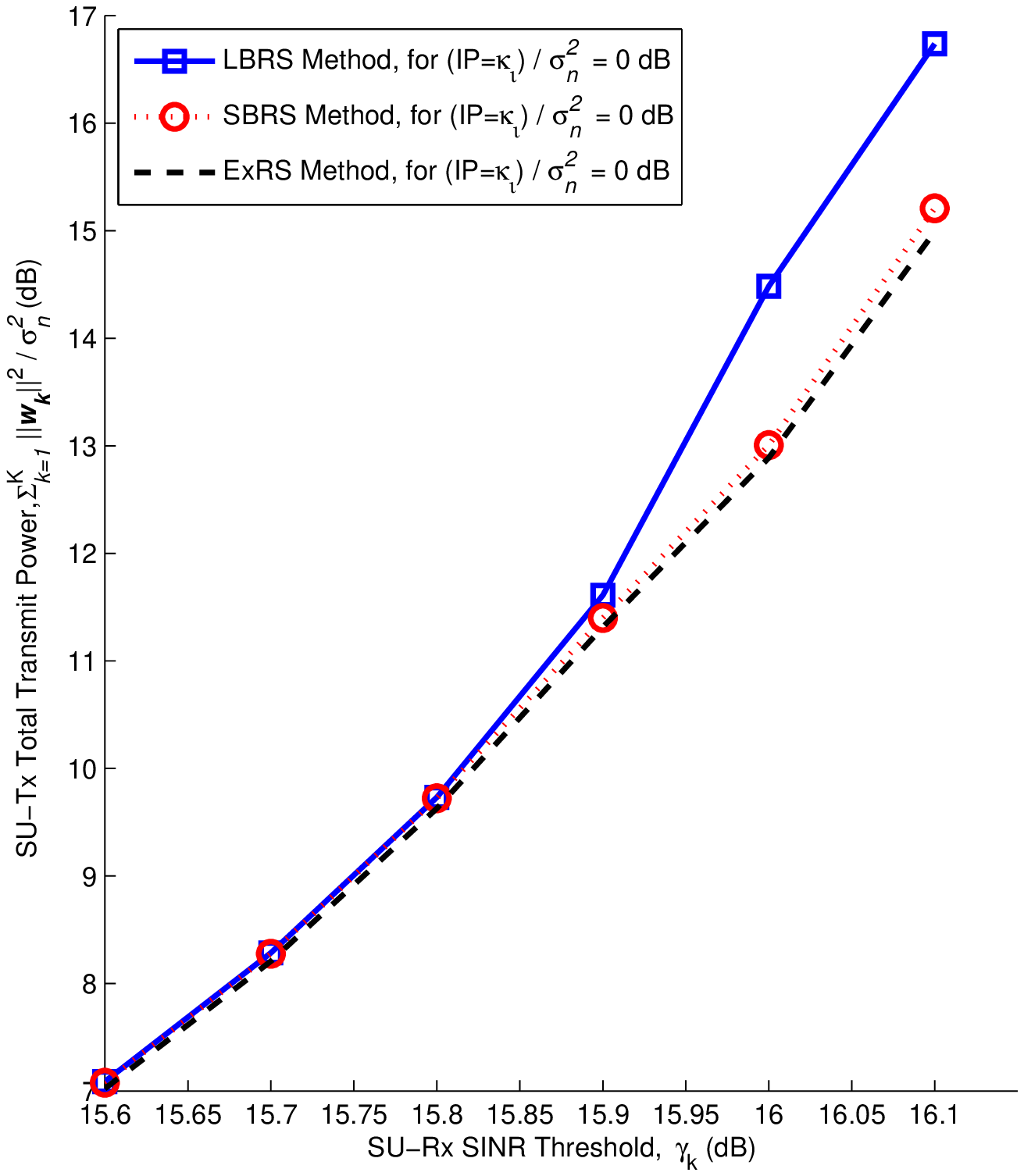}	\vspace{-0.8cm}\\
		\scriptsize (b) Relative IP Level of 0 dB, in High SNR Regime
	\end{tabular}
	\caption{The Total Tx Power vs. SINR Thresholds}
	\label{fig:TxP_vs_SINR}
\end{figure}

Also in Fig.~\ref{fig:TxP_vs_SINR} we have plotted the normalized total transmit power versus the SINR thresholds for different amounts of allowed normalized IP.
The normalized total transmit power is the ratio of total transmit power to the noise power and the normalized IP is defined using the same manner.
Both quantities are dimension-less and for better clarity are displayed in dBs.
As expected, ExCS is better than the other two methods.
In Fig.~\ref{fig:TxP_vs_SINR}-a, it is clear that for a relative IP level of $-$4 dB, ExCS transmits the lowest amount of power while SBCS requires to transmit a modest amount of power relative to LBCS, and finally LBCS requires to transmit the largest amount of power.
In this figure, it is also observed that for a relative IP level of 0 dB, the proposed ExCS is the best scheme to use to transmit power.
In Fig.~\ref{fig:TxP_vs_SINR}-b, we have plotted the same graph but in the higher SINR values.
In this range of SINR thresholds, all of the optimization problems with a relative IP level of $-$4 dB would be infeasible, so the graph is only provided for the relative IP level of 0 dB.
It is clear that in such scenarios, ExCS performs the best, although it should be noted that the performance of SBCS and ExCS are very close to each other.

\section{Conclusion\label{conclusion}}
The problem of robust downlink beamforming design in multiuser MISO cognitive radio networks is studied.
Particularly, a set up of $K$ SU-Rx's and $L$ PU-Rx's, all equipped with a single antenna is considered, and the SU-Tx has $N$ transmit antennas.
It is assumed that the relevant CSI is not known perfectly for both sets of users.
The uncertainty in the CSI is modeled using an Euclidean ball notation.
Three different approaches, namely LBCS, SBCS and ExCS, are presented which can be implemented efficiently.
The first solution is a SDP, while the later two solutions are the convex optimization problems.
Various simulation results are presented to evaluate the robustness of proposed methods.

\appendices
\section{Proof of Proposition 1 \label{app1}}
The Lagrangian function, using an arbitrary positive multiplier, $\lambda \geq 0$, is
\begin{align}
	L(\Mat{\Delta}_k,\lambda) &= \Tr{(\Mat{\tilde{H}}_k+\Mat{\Delta}_k) \Mat{W}_k} + \lambda (\| \Mat{\Delta}_k \|^2 - \epsilon_k^2) \nonumber \\
	&= \Tr{(\Mat{\tilde{H}}_k + \Mat{\Delta}_k) \Mat{W}_k} + \lambda (\Tr{\Mat{\Delta}_k \Mat{\Delta}_k^\dag} - \epsilon_k^2) .
\end{align}
By differentiating \cite{gesbert2007} of this Lagrangian function with respect to $\Mat{\Delta}_k^*$ and equating it to zero we will find the optimizer $\Mat{\Delta}_k$,
\begin{equation}
	\nabla_{\Mat{\Delta}_k^*} L(\Mat{\Delta}_k, \lambda) = \Mat{W}_k^\dag + \lambda \Mat{\Delta}_k = 0 ,
\end{equation}
which is annotated by $\Mat{\Delta}_k^{opt}$,
\begin{equation}
	\Mat{\Delta}_k^{opt} = -\frac{1}{\lambda} \Mat{W}_k^\dag.
\end{equation}
To eliminate the role of arbitrary parameter of $\lambda$, again, we differentiate the Lagrangian function with respect to this unknown parameter and then equate it to zero
\begin{equation}
	\nabla_{\lambda} L(\Mat{\Delta}_k, \lambda) = 0 ,
\end{equation}
to get the optimizer $\lambda$, annotated as $\lambda^{opt}$,
\begin{equation}
	\lambda^{opt} = \frac{1}{\epsilon_k} \| \Mat{W}_k^\dag \| .
\end{equation}
By combining these results, finally, we come up with
\begin{equation}
	\Mat{\Delta}_k^{opt} = - \epsilon_k \frac{\Mat{W}_k^\dag}{\|\Mat{W}_k\|} .
\end{equation}
To test if this solution belongs to a minimum, we should observe that the second derivative at the optimizer points should have a non-negative value:
\begin{equation}
	\nabla^2_{\Mat{\Delta}_k^*} L(\Mat{\Delta}_k^{opt},\lambda^{opt}) = \lambda^{opt} \geq 0 .
\end{equation}

To find the maximum of such a term, again, using a positive arbitrary Lagrangian multiplier, we build a Lagrangian function.
\begin{align}
	L(\Mat{\Delta}_k,\lambda) &= \Tr{(\Mat{\tilde{H}}_k+\Mat{\Delta}_k) \Mat{W}_i} - \lambda (\| \Mat{\Delta}_k \|^2 - \epsilon_k^2) \nonumber \\
	&= \Tr{(\Mat{\tilde{H}}_k + \Mat{\Delta}_k) \Mat{W}_i} - \lambda (\Tr{\Mat{\Delta}_k \Mat{\Delta}_k^\dag} - \epsilon_k^2) .
\end{align}
By differentiating it with respect to $\Mat{\Delta}_k$ and equating it to zero
\begin{equation}
	\nabla_{\Mat{\Delta}_k^*} L(\Mat{\Delta}_k, \lambda) = \Mat{W}_i^\dag - \lambda \Mat{\Delta}_k = 0  ,
\end{equation}
we will get
\begin{equation}
	\Mat{\Delta}_k^{opt} = \frac{1}{\lambda} \Mat{W}_i^\dag .
\end{equation}
Again, by differentiating the Lagrangian function with respect to $\lambda$ and equating it to zero,
\begin{equation}
	\nabla_{\lambda} L(\Mat{\Delta}_k, \lambda) = 0 ,
\end{equation}
we are able to get the optimizer.
\begin{align}
	\lambda^{opt} &= \frac{1}{\epsilon_k} \| \Mat{W}_i \| , \\
	\Mat{\Delta}_k^{opt} &= \epsilon_k \frac{\Mat{W}_i^\dag }{\|\Mat{W}_i \|} .
\end{align}
To prove if this solution belongs to a maximum, we should observe that:
\begin{equation}
	\nabla^2_{\Mat{\Delta}_k^*} L(\Mat{\Delta}_k^{opt},\lambda^{opt}) = -\lambda^{opt} \leq 0 .
\end{equation}

\section{Proof of Proposition 2 \label{app2}}
The Lagrangian multiplier is adopted again:
\begin{align}
	L(\Mat{\Delta}_k,\lambda) & = \Tr{(\Mat{\tilde{H}}_k + \Mat{\Delta}_k) \Mat{W}_k} - \gamma_k \sum_{\substack{i=1\\i\neq k}}^{K} \Tr{(\Mat{\tilde{H}}_k+\Mat{\Delta}_k) \Mat{W}_i} + \lambda (\Tr{\Mat{\Delta}_k \Mat{\Delta}_k^\dag} - \epsilon_k^2).
\end{align}
By differentiating this function and equating it with zero,
\begin{align}
	\nabla_{\Mat{\Delta}_k^*} L(\Mat{\Delta}_k,\lambda) = \Mat{W}_k^\dag - \gamma_k \sum_{\substack{i=1\\ i\neq k}}^{K} \Mat{W}_i^\dag + \lambda \Mat{\Delta}_k = 0 ,
\end{align}
we will come up with
\begin{align}
	\Mat{\Delta}_k^{opt} = -\frac{1}{\lambda} \left( \Mat{W}_k - \gamma_k \sum_{\substack{i=1\\ i\neq k}}^{K} \Mat{W}_i \right)^\dag ,
\end{align}
and to eliminate the $\lambda$,
\begin{align}
	\nabla_{\lambda} L(\Mat{\Delta}_k,\lambda) = \|\Mat{\Delta}_k\| - \epsilon_k = 0 ,
\end{align}
we will get
\begin{align}
	\lambda^{opt}&=\frac{1}{\epsilon_k} \|\Mat{W}_k - \gamma_k \sum_{\substack{i=1\\ i\neq k}}^{K} \Mat{W}_i\| ,\\
	\Mat{\Delta}_k^{opt} &= -\epsilon_k \frac{\left(\Mat{W}_k - \gamma_k \sum_{\substack{i=1\\ i\neq k}}^{K} \Mat{W}_i\right)^\dag}{\|\Mat{W}_k - \gamma_k \sum_{\substack{i=1\\ i\neq k}}^{K} \Mat{W}_i\|}  .
\end{align}
The second order differential test to prove that this solution belongs to a minimum, in this case, is also straight forward and is not included here.

%\section*{Acknowledgment}

\end{document}